\documentclass[12pt]{amsart}

\usepackage{amscd}
\usepackage{amsmath,amsfonts,amsthm,epsfig,latexsym,graphicx,amssymb,psfrag}
\usepackage[english]{babel}
\newtheorem{coro}{Corollary}
\newtheorem{defi}{Definition}[section]
\newtheorem{prop}{Proposition}[section]
\newtheorem{theo}{Theorem}
\newtheorem{lemm}{Lemma}[section]

\newtheorem*{theor}{Theorem}

\newtheorem{rem}{Remark}[section]

\def\R{I\kern -0.37 em R}
\def\N{I\kern -0.37 em N}
\def\Z{I\kern -0.37 em Z}

\def\supess_#1{\mathop{\rm supess}\limits_{#1}}
\def\infess_#1{\mathop{\rm infess}\limits_{#1}}

 \def\NN{{\mathbb N}} 

\def\QQ{{\mathbb Q}} \def\RR{{\mathbb R}} 
\def\TT{{\mathbb T}}

 \def\ZZ{{\mathbb Z}}

\begin{document}

\title[BS-actions on surfaces.]{\bf Actions of   Baumslag-Solitar groups   on
  surfaces.}
\author{Nancy Guelman and Isabelle Liousse} \thanks{This paper was partially
supported by Universit\'{e} de Lille 1, PEDECIBA, Universidad de la
Rep\'{u}blica and the PREMER project.}

\address{{\bf  Nancy Guelman}
IMERL, Facultad de Ingenier\'{\i}a, Universidad de la Rep\'ublica,
C.C. 30, Montevideo, Uruguay.  \emph{nguelman@fing.edu.uy}, }

\address{{\bf    Isabelle Liousse}, UMR CNRS 8524, Universit\'{e} de Lille1,
59655 Villeneuve d'Ascq C\'{e}dex,   France.  \emph {liousse@math.univ-lille1.fr},  }

\begin{abstract} Let $BS(1, n) =< a, b \ | \ aba^{-1} = b^n >$  be the solvable
  Baumslag-Solitar group, where $ n\geq 2$. It is known that $BS(1, n)$
  is isomorphic to the group generated by the two affine maps  of the real line:
  $f_0(x) = x + 1$ and $h_0(x) = nx $.

This paper deals with the dynamics of actions of  $BS(1, n)$
on closed orientable surfaces. We exhibit  a smooth $BS(1,n)$ action without
finite orbits  on $\TT ^2$,
we study  the dynamical behavior  of it  and of its $C^1$-pertubations and we prove that it  is not locally rigid.

We develop a general dynamical study for faithful  topological
$BS(1,n)$-actions on closed surfaces $S$. We prove that such
actions $<f,h \ \vert \  h \circ f \circ h^{-1} = f^n>$ admit a minimal set included in $fix(f)$,  the set  of
fixed points of $f$, provided that  $fix(f)$ is not empty.

When $S= \TT^2$, we show that there exists a positive  integer  $N$, such that
$fix(f^N)$ is non-empty and contains a minimal set of the action. As
a corollary, we get that there are no minimal faithful topological
actions of $BS(1,n)$ on  $\TT^2$.

When the surface $S$ has genus at least 2,  is closed and orientable,   and $f$ is isotopic to identity, then  $fix(f)$ is non empty and contains a  minimal set of the action. Moreover if the action  is $C^1$  then $fix(f)$  contains any minimal set.

\end{abstract}

\maketitle
\date{}
\section {Introduction and statements}

\medskip

An important question on  group actions is existence and stability of global fixed points. For Lie group actions,
it was shown by Lima \cite{Li}  that any action of the abelian Lie group $\RR ^n$ on a surface with
non-zero Euler characteristic  has a global fixed point. This result was later extended
by Plante \cite{P} to nilpotent Lie groups.  On the other hand, Lima \cite{Li} and  Plante\cite{P} proved that
the solvable Lie group $GA(1,\RR)$ acts without fixed points on every compact surface.

For discrete group actions, Bonatti \cite{B}  showed that any $\ZZ^n$ action on surfaces with non-zero Euler characteristic generated
by diffeomorphisms $C^1$ close to the identity has a global fixed point.  Druck, Fang and Firmo \cite{DFF}  proved a discrete version of Plante's theorem.

This paper deals with the dynamics of  actions of the solvable
Baumslag-Solitar group, $BS(1, n) =< a, b \ | \ aba^{-1} = b^n >$,
where $ n\geq 2$, on closed surfaces.

\medskip

It is well known that  $BS(1,n)$  has many  actions on $\mathbb R$.
 The standard  action on $\mathbb R$ is
 the action   generated by the two affine maps $f_0(x) = x + 1$ and $h_0(x) = nx $ (where $f_0 \equiv b$ and $h_0 \equiv a $).

Actions of solvable groups on one-manifolds have been studied by  Plante \cite{P}, Ghys \cite{Gh}, Navas \cite{Na}, Farb and Franks\cite{FF},
 Moriyama\cite{Mo} and Rebelo and Silva\cite{RS}.  In \cite{FF}, as a corollary
of a result of M. Shub \cite{Sh} on expanding maps, they showed the following:

\begin{theor}{\text {(Farb-Franks-Shub).}}

There are neighborhoods of $f_0$  and   $ h_0$ in the uniform $C^1$ topology
such that whenever $f$ and $h$ are chosen from these respective neighborhoods
and  the group generated by $\{f,h\}$ is isomorphic to $ BS(1, n)$  then the
perturbed action is topologically conjugate to the original action.

\end{theor}

In contrast, Hirsch (see \cite{Hi})  has found analytic
actions of $BS(1, n)$  on $\mathbb R$  which are not topologically conjugate to
the standard action (but they are semiconjugate).

\medskip

\begin{defi}

The {\bf standard $ BS(1,n)$-action on $S^1= \RR \cup {\infty}$} is the action  generated by the two Moebius maps
 $f_0(x) = x + 1$ and $h_0(x) = nx $. It has a global fixed point at ${\infty}$.
\end{defi}


Farb-Franks-Shub Theorem remains true for the standard $
BS(1,n)$-action on $S^1$, since a faithful BS-action $C^1$ close to
the standard action on $S^1$ always has a global fixed point (the proof of this fact is analoguous to the proof of lemma \ref{pers}) and therefore it can be seen as an action  on $\mathbb R$.
\bigskip

More recently, L.  Burslem and A. Wilkinson \cite{BW} gave a classification  (up to
conjugacy) of real analytic actions of $BS(1, n)$ on $S^1$. In particular,  they
 proved that every representation  of $BS(1, n)$ into $ {\text {Diff}}^\omega (S^1)$
is $ C^\infty$-locally rigid and that for   each $r\geq 3$  there are
analytic actions of $BS(1, n)$  that are $C^r$-locally rigid, but not
$C^{r-1}$-locally rigid (for the definition of local rigidity see section
\ref{preli}).

This results are proved by using a dynamical approach. The dynamics
of $C^2$ $BS(1,n)$-actions on $S^1$ is now well understood, due to
Navas work on solvable groups of circle diffeomorphisms (see
\cite{Na}).
 In particular, Burslem and Wilkinson (see \cite{BW}) proved that any $C^2$ faithful $BS(1,n)$-action on $S^1$ admits a finite
 orbit. Recently we have extended this result to $C^1$ case (see
 \cite{GL}). Also, we have proved that any $C^1$ faithful  $BS(1,n)$-action on
 $S^1$ is semiconjugated (up to passing to a finite index subgroup) to the standard one.

\bigskip

The dynamical situation of $BS(1,n)$-actions  on closed surface (even on
$\TT^2$)   is more complicated.  Non trivial examples of $BS(1,n)$ actions on
closed surfaces  can be constructed  using actions of the affine real group
 $GA(1,\RR):= \{ x\mapsto \alpha x + \beta, \ \   \alpha, \beta\in \RR, \
 \alpha>0 \}$. Actions of $GA(1,\RR)$ on closed surfaces have been studied by Plante-Thurston \cite{PT}, Plante \cite{P}, Belliart-Liousse      \cite{BL}.

On the other  hand, $C^1$ faithful $BS(1,n)$-action on $\TT^2$  can
be constructed by using products of smooth  actions on $S^1$. In the
case where the circle actions are both  $C^1$ and faithful,  finite
orbits always exist. But,   faithful smooth $BS(1,n)$-actions  on
$\TT^2$ can be obtained  as the product of a faithful
$BS(1,n)$-action and a non faithful one. In this case, finite orbits
may not exist.

\smallskip

An important family of such examples is given by $<f_0, h_k>$, where :

$$f_0(x,\theta)=(x+1,\theta)  \mbox{ and } h_k(x,\theta)=( nx, k(\theta) ),$$
 where $x\in \RR \cup \infty$,  $\theta \in S^1$ and $k$ is any (orientation preserving) circle
 homeomorphism.
\smallskip

\medskip

In section 3, we explain the construction of these examples  and exhibit  a  faithful smooth action  of  $BS(1,n)$ on $\TT^2$
without finite orbits that can be considered  as ``the standard $BS$-action'' on $\TT^2$. More precisely,

\begin{defi}  The {\bf standard $BS(1,n)$-action on $\TT^2$} is the action generated by :
 $$f_0(x,\theta)=(x+1,\theta)$$ and $$h_0(x,\theta)=( nx, \ln(n)+\theta),$$
 where $x\in \RR \cup \infty$ and $\theta \in S^1$.

\end{defi}

Our first result is the following:

 \begin{theo} \

 The group $<f_0,h_k>$ generated by $f_0$ and $h_k$ is isomorphic to $BS(1,n)$.

If the rotation number of $k$ is rational, there exist finite BS-orbits.

If the rotation number of $k$ is irrational, there are no finite BS-orbits and the unique minimal set for the BS-action is included in $\infty \times S^1=fix(f_0)$.

\label{f0h0}
\end{theo}

\begin{coro}\label{C2}
There exist $C^{\infty}$  faithful BS-actions arbitrary
$C^{\infty}$-close to  the standard torus BS-action  $<f_0,h_0>$
that are not topologically  conjugate to $<f_0,h_0>$.
 \end{coro}

This implies  that the  standard BS-action on $\TT^2$ does not
satisfy the rigidity properties described in the Farb-Franks-Shub
theorem for the standard BS-action on $S^1$.

This property can also be compared to the rigidity result recently
proved by Mc Carthy~: {\it ``The trivial $BS(1,n)$-action on a
compact manifold  does not admit  $C^1$ faithful perturbations"}
(see \cite{McC}).

\bigskip

Then we consider  perturbed actions  of the standard  one. In
particular, we prove that there exists either a finite orbit or a unique minimal set. Recall that  a {\bf  minimal set} for an action of a group $G$ on a  compact metric space $X$ is a non-empty
closed $G$-invariant subset of $X$ such that if $K\subset M$ is a closed $G$-invariant set then either $K=M$ or $K=\emptyset$.

\bigskip

Let $\mathcal C_1$ and $\mathcal C_2$ be the circles defined by $\mathcal C_1=\infty \times S^1$ and $\mathcal C_2=0\times S^1$.
Note that both circles are $h_0$-invariant.

\begin{theo}\label{minset}

Let us consider a BS-action $<f,h>$ on $\TT ^2$ generated by  $f$
and $h$ sufficiently $C^1$-close to  $f_0$ and $h_0$ respectively.
Then:

\begin{enumerate}

\item  there exists two  circles $\mathcal C'_1$ and $\mathcal C'_2$
close to $\mathcal C_1$ and $\mathcal C_2$ respectively which are
$h$-invariant. Moreover, the $\omega _h$-limit set of any point in
$\TT^2\setminus\mathcal C'_2 $ is included in  $\mathcal C'_1$ and
the $\alpha _h$-limit set of any point in $\TT^2\setminus\mathcal
C'_1$ is included $\mathcal C'_2$.

\item  the set of $f$-fixed points is not empty and it is contained in the
 circle  $\mathcal C'_1$.

\item    either :

\begin{enumerate}

\item  there  exist finite BS-orbits contained in  $\mathcal C'_1$, or

\item the action has a unique minimal set $M$ which is included in $\mathcal C'_1$ (and in the set of $f$-fixed points). Moreover, $M$ is  either $\mathcal C'_1$ or a Cantor set.

\end{enumerate}

\end{enumerate}

\end{theo}

\bigskip

We check that the ``standard  action'' on $\TT^2$ satisfies   item
3(b) but  in the proof of Corollary \ref{C2} we exhibit  $C^{\infty}$-perturbations of it that    have
a different dynamical behavior: they satisfy item 3(a). In section 6, we  exhibit an example of an action with a $C^1$ persistent global fixed point. More precisely, we construct an action with  fixed point satisfying that any  $C^1$-perturbation of it  also has fixed point.

\medskip

On the other hand,  we develop a general dynamical study for  faithful
$BS(1,n)$-actions on closed surfaces. From now on, let us consider $f$ and $h$ two
homeomorphisms that generate a   $BS(1,n)$-action,  that is, $h\circ f\circ h ^{-1} = f^n$.

Our first ``dynamical'' result on the torus concerns the rotation set of $f$ (for the
definition  see Section \ref{preli}).

 \begin{theo}\label{rot}
 Let  $<f,h>$ be a faithful action  of $BS(1,n)$ on   $\TT^2$. Then  there exists a
  positive  integer  $N$, such that  $f^N$  is isotopic to identity and has a lift whose rotation
   set is the single point $\{(0,0)\}$. Moreover, the set of $f^N$-fixed points denoted by  $fix(f^N)$ is non-empty.

   \end{theo}

\begin{rem} In section 3, we exhibit two diffeomorphisms $F$ and $H$ generating a faithful action  of $BS(1,n)$ on   $\TT^2$,
where  $F$ admits periodic orbits but it does not have fixed points.
\end{rem}
\medskip

Since the group $<f^N,h>$ is isomorphic to $ BS(1,nN)$, Theorem \ref{rot}
allows us to restrict our study on the torus to the case where $f$
is isotopic to identity, the rotation set of a lift of $f$ is
$\{(0,0)\}$ and $f$ has fixed points.  In this situation we prove
that there exists a $BS$-minimal set included in the set of
$f$-fixed points. More precisely, we prove the more general
following statement.

\medskip

\begin{theo}\label{minBS}

Let $X$ be a compact metric space and $<f,h>$ be a representation of
$BS(1,n)$ in $ Homeo(X)$.

\noindent (a) If $fix (f)$  is  non-empty,  then:

\begin{enumerate}

\item  If $x_0\in fix (f)$ then   $\alpha_h(x_0)$   is contained in
$fix(f)$.

\item There exists an $BS$-minimal set  included in $fix(f)$. Moreover, this $BS$-minimal set coincides with a $h$-minimal set in $fix(f)$.

\item If the set of $f$-fixed points is finite then the action admits a global
 finite orbit.

\item Let $\mathcal{M}$ be an $BS$-minimal set satisfying $\mathcal{M} \cap
  fix(f)\neq \emptyset$, then  $\mathcal{M} \subset fix(f)$.


\end{enumerate}

\medskip

\noindent (b) If the set  of periodic points of $f$, $per(f)$,  is
non-empty, then there exist a positive integer $N$ and a
$BS$-minimal set, $M$,  such that $M \subset fix f^N$.

\end{theo}

As a consequence of item (b) of Theorem \ref{minBS}, Theorem
\ref{rot} and the fact that $<f,h>$ is a faithful representation of
$BS(1,n)$, we have the following:

\begin{coro}\label{nomin}

Let $X$ be a compact metric space and $<f,h>$ be a faithful representation of $BS(1,n)$ in $Homeo(X)$ such that
$Per (f)$  is  non-empty. Then the action of $<f,h>$ is not  minimal.

\noindent In particular:

\begin{enumerate}
\item There is no faithful minimal action of  $BS(1,n)$ by homeomorphisms on $\TT ^2$.
\item Let $\Sigma$ be a compact surface of non zero Euler characteristic.  A faithful  topological $BS(1,n)$-action  $<f,h>$ on  $\Sigma$ is not minimal, provided that  $f$ is isotopic to identity.
\end{enumerate}

\end{coro}
\begin{rem} Item (2) is a consequence of Theorem \ref{minBS} (a)(2) and Lefschetz's fixed point theorem. \end{rem}

\medskip

When $<f,h>$ is a topological action of $BS(1,n)$ on a closed
surface  $S$  satisfying that  any $f$-invariant probability has support
included in the set of $f$-fixed points, we prove the following:

\begin{theo}
\label{mu}
Let   $S$ be a closed  orientable surface and $<f,h>$  be a representation of  $BS(1,n)$ in  $Homeo(S)$.
Suppose that for any $f$-invariant probability measure $\mu$, $supp(\mu) \subset fix(f)$. Then:

\begin{enumerate}
\item Any $f$-minimal set is a fixed point. The set  of periodic points of $f$, $per(f)$, coincides with the set $fix(f)$.

\item Any $BS$-minimal set is  included in $fix(f)$. Moreover, any $BS$-minimal set coincides with a $h$-minimal set in $fix(f)$.

\item  Topological entropy of $f$, $ ent_{top}(f)=0$.
 \end{enumerate}

\end{theo}

For next corollary, that we prove using Theorem 1.3 of
\cite{FH}, we need the following:
\begin{defi}
Let $g\in Diff^1(S)$, an $N$-periodic point $x_0$ is called {\bf elliptic}
if the eigenvalues of the differential of $g$ at $x_0$, $Dg^N(x_0)$,
have module $1$.
\end{defi}

\begin{coro}
\label{C1}
Let  $S$ be a closed  orientable surface and   $<f,h>$  be a representation of  $BS(1,n)$ in  $Diff^1(S)$  such that:

\begin{itemize}
\item If $S$ has genus at least 1, $f$ is isotopic to identity.

\item  If $S=S^2$, some iterate of $ f$ has at least three fixed points.
 \end{itemize}

Then there exists a positive integer $N$ such that :
\begin{enumerate}
\item Any $f$-minimal set is a periodic point. The set  of periodic points of $f$, $per(f)$, coincides with the set $fix(f^N)$.

\item Any $BS$-minimal set is  included in $fix(f^N)$. In fact, any $BS$-minimal set is  included in a subset of $f$-elliptic points in $fix(f^N)$.

\item  $ ent_{top}(f)=0$.

\end{enumerate}
In addition, if $S$ has genus at least 2, $N=1$.
\end{coro}


Finally, we have the following open questions:

\begin{enumerate}
\item Does it
exist  a faithful  action of $<f,h>=BS(1,n)$  on $\TT^2$ with $h$ non isotopic to identity? We know that there does not exist representation of $BS(1,n)$  into  $Aff(\TT^2)=SL(2,\mathbb Z)\ltimes \mathbb R^2$ the group consisting of maps :$g(x,y)= A.(x,y) + V$, where $A\in SL(2,\mathbb Z)$ and $ V\in \mathbb R^2$.

\item Does it exist a faithful continuous action of $<f,h>=BS(1,n)$  on $\TT^2$ with minimal sets outside $per(f)$ ?

\item Is the product action on  $ (\RR \cup \infty ) \times S^1$  generated by $f_0(x,\theta)=(x+1,\theta)$
  and $h_0(x,\theta)=( nx, k(\theta) ),$ where $k$ is a circle  $north$-$south$ diffeomorphism topologically rigid ?


\end{enumerate}

In Section 2, we give definitions, properties and basic tools that
we use in the rest of the paper. We exhibit examples of $BS(1,n)$
acting on $\TT^2$, and Theorem 1 and Corollary \ref{C2} are proved
in Section 3. The goal in Section 4 is proving Theorem 3. In Section
5 we prove Theorems \ref{minBS}, \ref{mu} and Corollaries
\ref{nomin} and \ref{C1}. In Section 6,  we consider perturbations
of the standard $BS(1,n)$-action on $\TT^2$: we describe their
minimal sets by proving Theorem \ref{minset}. We also construct an action with a persistent global fixed point.

\section{Definitions-Notations}\label{preli}

\subsection{Isotopy class of torus homeomorphims.}\

\medskip

We denote by $Homeo_{\mathbb Z ^2}(\RR ^2)$ the set of homeomorphisms $F : \RR
^2 \rightarrow \RR ^2$ such that $F(\mathbb{Z }^2 )\subseteq \mathbb Z ^2$ and $Homeo^0_{\mathbb Z ^2}(\RR ^2)$ the set of homeomorphisms $F : \RR
^2 \rightarrow \RR ^2$ such that $F(\mathbb{Z }^2 )\subseteq \mathbb Z ^2$ and $F(x + P) = F(x) +P$, for all $x\in \RR^2$
and $P\in \ZZ ^2$.

\smallskip

\noindent Note that, a lift  of a 2-torus homeomorphism
isotopic to identity  belongs to $Homeo^0_{\mathbb Z ^2}(\RR ^2)$. Conversely,
if a  2-torus homeomorphism admits a lift $F\in Homeo^0_{\mathbb Z ^2}(\RR ^2)$
then it is isotopic to identity.

Let $ g : \TT^2 \rightarrow \TT^2$ be a homeomorphism and let $G : \RR^2
\rightarrow \RR^2$ be a lift of $g$.  We can associate to $G$ a linear map
$A_G $  defined by :

$$  G( p + (m,n)) = G( p )  + A_G (m,n), {\text { for any }}  m,n  {\text
  { integers.} }$$

By definition, it is clear that  $G\in Homeo^0_{\ZZ ^2} (��\RR ^2)$ (that is
$g$ is isotopic to identity) if and only if  $A_G= Id$.

This map satisfy the following properties :

\begin{enumerate}

\item  $A_G$  does not depend neither on the integers $m$ and $n$ nor on the lift
$G$ of $g$.  In fact,  $A_G$ is the morphism  induced by $g$ on the first
homology group of $\TT ^2$.  So we can also denote $A_g$ for $A_G$ and we will
use both notations.

\item  $A_{G\circ F} = A_{G} \circ A_{ F}$ and  $A_{G^{-1}} = A_{G} ^{-1}$,

\item  $A_G \in GL(2,\ZZ)$,  in particular

\item  $det A_G = +1 $ or $-1$.

\end{enumerate}

\medskip

\subsection{Rotation set and  rotation vectors.}

\subsubsection{\bf Definitions} \

Let  $f$ be a  2-torus homeomorphism
isotopic to identity. We denote by  $\tilde f$ a lift of to  $\mathbb R ^2$. We call  {\bf ${\tilde
    f}$-rotation set} the subset of $\mathbb R ^2$ defined by

$$ \rho(\tilde f)= \bigcap_{i=1}^{\infty}\overline{\bigcup_{n\geq i} \left\{\frac {
    \tilde f   ^{n} (\widetilde{x}) - \widetilde{x}} {n}, \ \  \widetilde{x} \in \mathbb{R}^2\right\}}.$$

Equivalently, $(a,b) \in \rho(\tilde f)$ if and only if there exist
sequences $(\widetilde{x_i})$ with $\widetilde{x_i} \in
\mathbb{R}^2$ and $ n_i \to \infty$ such that
 $$(a,b)= \lim_{i\to \infty} \frac { \tilde f
  ^{n_i} (\widetilde{x_i}) - \widetilde{x_i}} {n_i}.$$

\medskip

Let $\widetilde x$  be in $\mathbb R^2$. The  {\bf $\widetilde{f}$-rotation
  vector   of $\widetilde x$} is the
2-vector defined by $\displaystyle \rho(\widetilde{f},\widetilde x)=\lim_{n\to
\infty} \frac { \tilde f
  ^{n} (\widetilde{x}) - \widetilde{x}} {n} \in \mathbb{R}^2$  if this limit exists.

  From now on, we use  both $\tilde f$ or $F$ for a lift of $f$ to  $\mathbb R ^2$.

\bigskip

\subsubsection {\bf Some classical properties and results on  the rotation
set.}

Let $f$ be a  2-torus homeomorphism  isotopic to the identity and
 $\widetilde{f}$ be a lift of $f$ to $\RR^2$.

\medskip
\begin{itemize}
 \item Let $\widetilde x\in \mathbb R^2$ be such that
   $\rho(\widetilde{f},\widetilde x)$ exists. Then
   \ \ $\rho(\widetilde{f},\widetilde x)\in \rho(\widetilde{f})$.

\item If $\widetilde{f}$ has a  fixed point   then  $(0,0) \in
\rho(\widetilde{f})$.

\item Misiurewicz and Ziemian (see  \cite{MZ89}) have  proved that:
\begin{quote}
\begin{enumerate}

\item $\rho(\tilde f^n ) = n \rho(\tilde f)$

\item  $\rho(\tilde f +(p,q)) = \rho(\tilde f ) +(p,q) $,

\item the rotation  set is a  compact convex subset of
$\mathbb R^2$.

\end{enumerate}
\end{quote}
\end{itemize}


\subsection{ \bf  $ C^r$-local rigidity, where $r \in \NN \cup \{ \infty, \omega\}$.}

%
%
%
%
%
%
%
%
%
\begin{defi}
An action $<f_1, h_1>$ of $BS(1,n)$ on a smooth manifold  is {\bf $C^r$-locally rigid}
($r\in \mathbb{N}\cup \infty\cup \omega$) if there
  exist  neighborhoods of  $f_1$ and $h_1$ in the $C^1$-topology such that whenever $f$ and $h$ are $C^r$ maps chosen from these neighborhoods  and the group generated by $<f, h>$ is  isomorphic to $BS(1,n)$, then the perturbed action
    is $C^r$ conjugate to the original one, that is there exists a $C^r$-diffeomorphism $H$ such that $ H\circ f \circ H^{-1}=f_1$ and $ H\circ h \circ H^{-1}=h_1$.

\end{defi}













\subsection{\bf  $BS(1,n)$-actions} {Consequence of the conjugation between $f^ n$ and $f$.}

As consequences  of the group-relation  $h\circ f\circ h^{-1}=f^n$, we get
  easily the
following two propositions :

\medskip

\begin{prop} \ Let $f$ and $h$ be homeomorphisms satisfying $h\circ f\circ h^{-1}=f^n$, then

\begin{enumerate}

\item  $h \circ f^p \circ h^{-1}= f^{np} $, for all integer $p$,
\item  $h^p \circ f \circ h^{-p} = f ^{n^p} $, for all positive integer $p$.

\end{enumerate}

\end{prop}

\medskip

\begin{prop}\

Let $f$ and $h$ be as in the previous proposition, then
\begin{enumerate}

\item $h(fix (f)) = fix (f^n)$,

\item  Let $per(f)$ be the set of periodic points of $f$, then  $h(per(f)) =
  per (f^n)$.  More precisely,  if  $x$ is an  $f^p$ fixed point then
$h(x)$ is an $(f^n)^p=f^{np}$ fixed point.

\item  If $M_f$ is an $f$-minimal set then $h(M_f)$ is a minimal set of $f^n$.

\item  Let $ent(f)$ be the topological entropy of $f$. Then $ent(f)$ is $0$ or
$\infty$.

\end{enumerate}

\end{prop}

\medskip

\noindent {\bf Proof of (4).} Since $ent(f^n) =  n .  ent(f)$ and $ent(h\circ f \circ h ^{-1})=ent(f)$ the
possible values for $ent(f)$ are  $0$ or $\infty$.

\medskip

\section{Examples  of $BS(1,n)$-actions on $\TT ^2$}  \label{example}

In this section we will exhibit examples of $BS(1,n)$-actions on $\TT ^2$.

\subsection {Product of  faithful actions on $S^1$.} \

\medskip

Let $<f_i, h_i>$, $i=1,2$  be two $C ^1$ actions of $BS(1,n)$ on $S^1$, we construct
an action of $BS(1,n)$ on $\TT ^2$ by setting :
$f= (f_1, f_2)$ and $h=(h_1, h_2)$. Clearly, the $<f,h>$-orbit of a point
$x=(x_1,x_2)$  in $\TT ^2$ is the product  of the $<f_1,h_1>$-orbit of $x_1$
and the  $<f_2,h_2>$-orbit of $x_2$.

According to \cite{GL}, there exists a finite  $<f_i,h_i>$-orbit at
some point $y_i \in S^1$, hence the  $<f,h>$-orbit of the point
$y=(y_1,y_2)$  is finite.

\medskip
The following two sections show examples of $BS$-actions on $\TT^2$
without finite orbits.

\subsection {Product of  non faithful actions on $S^1$.} \

We construct faithful $BS(1,n)$-actions  without finite orbits as product of a faithful
circle action  and a non faithful one.

Let $<f_1, h_1>$  be a faithful action of $BS(1,n)$ on $S^1$ and $k$
be a circle homeomorphism. We construct a faithful  action of
$BS(1,n)$ on $\TT ^2$ by setting: $f= (f_1, Id)$ and $h=(h_1, k)$.
Clearly, if $k$ has no finite orbit, there is no global finite
orbit.

\medskip

\subsection {Actions that come from  actions of the affine group of the real line.}

\medskip

\subsubsection{\bf Actions of $GA(1,\mathbb R)$ and induced $BS(1,n)$-actions on the circle.} \

\medskip

Identifying the affine real map $x \mapsto ax+b$ with $(a,b)$,  the affine group of the real line, $GA(1,\mathbb R)$, is the group $\mathbb R_{>0} \times \mathbb R$ endowed with the product $(a,b)\times (a',b')= (aa', ab'+b)$.

The Baumslag-Solitar group $BS(1,n)$ can be seen as  the  subgroup generated by the elements $(1,1)$ and $(n,0)$.

Let   $\Phi :  GA(1,\mathbb R)  \rightarrow Diff^{r} (M)$ be an
action of $GA(1,\mathbb R)$,   the {\bf  induced  $BS(1,n)$-action}  is the
restriction of $\Phi$ to $<(1,1), (n,0)>$.

 \bigskip

{\bf The standard actions.}

\begin{defi}

The  {\bf standard action of $GA(1,\mathbb R)$ on the circle} is the action by Moebius maps on the projective line, that is :

$$\Phi^{stand} : \left \{ \begin{array}{ll} GA(1,\mathbb R)  &\rightarrow Diff^{\omega}_+ (S^1) \\
                                      (a,b) & \mapsto  \Phi^{stand}_{(a,b)} \end{array} \right . , $$

{\text { where }}

$$\displaystyle \Phi^{stand}_{(a,b)}:\left\{ \begin{array}{ll}  \mathbb R \cup \{\infty\}  &\rightarrow \mathbb R \cup \{\infty\} \\                                                                            x & \mapsto ax +b \end{array} \right . $$
\end{defi}

This action is faithful and has a global fixed point at $\infty$.

\bigskip
\begin{defi}

The  {\bf standard action of $BS(1,n)$ on the circle} is the induced $BS(1,n)$-action, it is generated by the two Moebius maps  $f_0(x)=\Phi^{stand}_{(1,1)} (x) = x+1$ and $h_0(x)=\Phi^{stand}_{(n,0)} (x) = nx$.
\end{defi}

It is faithful and has a global fixed point at $\infty$. Moreover
$f_0$ has a unique fixed point at $\infty$ that is elliptic and
$h_0$ has two hyperbolic fixed points : $\infty$ that is an
attractor and $0$ that is a repeller.

The orbit of a point $x$ is explicit : $\mathcal O (x) = \{n^k x + w, k\in \mathbb Z, w \in \mathbb Z[ \frac {1}{n}]\}$.

All orbits are dense except the orbit of the global fixed point $\infty$.

\begin{rem}
Applying the change of coordinate $x=\tan (\frac {u}{2})$ the standard $GA(1,\mathbb R)$-action is given by :

$$\displaystyle \Phi^{stand}_{(a,b)}:\left\{ \begin{array}{ll}  {[-\pi, \pi]}/ {(-\pi \sim \pi)}  &\rightarrow {[-\pi, \pi]} /{(-\pi \sim \pi)} \\  u &\mapsto 2\arctan (a\tan (\frac {u}{2}) +b ) \end{array} \right . $$
\end{rem}

\bigskip

{\bf  Non faithful actions.} \

A family of non faithful action is given by :

$$\Phi^{deg} : \left \{ \begin{array}{ll} GA(1,\mathbb R)  &\rightarrow Diff^{\omega}_+ (S^1) \\
                                      (a,b) & \mapsto  \varphi_{\ln a }\end{array} \right . , $$

 where $(\varphi_t)$ is any flow on the circle.

\bigskip

The induced $BS(1,n)$-actions are the actions generated by
$f(\theta) = \theta$ and $ h(\theta)= \varphi_{\ln n}(\theta)$.

\medskip

\begin{rem} There exist actions that do not come from actions of the affine group of the real line:
There exist (even orientation preserving) circle homeomorphisms
which do not embed in a continuous flow (see \cite{Zd}). However,
the family $ <f(\theta) = \theta, h(\theta)= \varphi_{\ln n}(\theta)>$, where $\varphi$ is a flow,   extends to  actions
$ <f(\theta) = \theta, h(\theta)=k(\theta)>$, where $k$ is any circle homeomorphism.

It is easy to see that $h\circ f\circ h^{-1} =f^n$ and that these
actions are not faithful and have the dynamics of $k$ : $\mathcal O
(x) = \{k^n(x), n\in \mathbb Z\}$.

\end{rem}

\subsubsection{\bf Actions of $GA(1,\mathbb R) $ and induced $BS(1,n)$  on the 2-torus.} \

 Taking the product of the standard action with a non faithful action of $GA(1,\mathbb R)$ on the circle, we get a family  of faithful $GA(1,\mathbb R)$-actions on the 2-torus :

$$\Phi^{\varphi} : \left \{ \begin{array}{ll} GA(1,\mathbb R)  &\rightarrow Diff^{\omega}_+ (\TT^2) \\
                                      (a,b) & \mapsto  \Phi^{\varphi}_{(a,b)} \end{array} \right . , $$

where

$$\displaystyle \Phi^{\varphi}_{(a,b)}:\left\{ \begin{array}{ll}  (\mathbb R \cup \{\infty\} )\times S^1 &\rightarrow (\mathbb R \cup \{\infty\})\times S^1 \\                                                                            (x,\theta) & \mapsto (ax +b , \varphi_{\ln a} (\theta) ) \end{array} \right . $$

and $(\varphi_t)$ is any flow on the circle.

\bigskip

The ``extended" induced $BS(1,n)$-actions are the actions generated
by $f_0(x, \theta) = (x+1, \theta)$ and $ h_k(x, \theta)= (nx,
k(\theta))$, where $k$ is any circle orientation preserving
homeomorphism. They are faithful since they are products of two
actions, and one of them is faithful. Their dynamics depend on $k$.

\bigskip

\begin{defi} \

The  {\bf standard action of $GA(1,\mathbb R)$ on the 2-torus} is
the action $\Phi^{\varphi} $, where  $\varphi_t (\theta) = \theta +
t$ is the  flow of the circle rotations, that is given by
$$\displaystyle \Phi^{stand}_{(a,b)}:\left\{ \begin{array}{ll}  (\mathbb R
\cup \{\infty\} )\times S^1 &\rightarrow (\mathbb R \cup
\{\infty\})\times S^1 \\
(x,\theta) & \mapsto (ax +b , \theta + {\ln a} ) \end{array} \right
. $$

The  {\bf  standard action of $BS(1,n)$ on the 2-torus} is the
induced action, that is the action generated by $f_0(x, \theta) =
(x+1, \theta)$ and $ h_0(x, \theta)= (nx, \theta +\ln n)$.

\end{defi}

This $GA(1,\mathbb R)$-action has no global fixed point and it has an
1-dimension circular orbit $\{\infty\} \times S^1$.

This $BS(1,n)$-action has no finite orbit, the restriction of $h_0$
to $\infty \times S^1$ is the irrational rotation by $\ln n$. The
unique minimal set is $\infty \times S^1$.

\subsection {Proof of Theorem \ref{f0h0} and Corollary \ref{C2} } \

\noindent {\bf Proof of Theorem \ref{f0h0}.}

We considered actions  on the 2-torus  generated by $f_0(x, \theta) = (x+1, \theta)$ and $ h_k(x, \theta)= (nx, k(\theta))$.

In the previous section, we have seen that  these actions are
faithful $BS(1,n)$-actions, since they are  products of two
$BS(1,n)$-actions on $S^1$, where one of them is faithful.

Note that the set of $f_0$-fixed points, $fix(f_0)$  is the circle $\mathcal C_1:=\infty \times S^1$ and any horizontal circle  $(\mathbb R \cup \{\infty\}) \times \theta_0$ is $f_0$-invariant. The  circles $\mathcal C_1=\infty  \times S^1$ and   $\mathcal C_2:=0
 \times S^1$ are $h_k$-invariant. The restriction of $h_k$ to these circles is the homeomorphism $k$.

\begin{itemize}

\item If the rotation number of $k$ is rational, then  there exists a point in  $\infty \times S^1$ with a $h_k$-finite orbit. As it is $f$-fixed, its   BS-orbit is finite.


\item If  the rotation number of $k$  is irrational (for example for the standard action), there are neither fixed points nor  periodic points of $h_k$.
Therefore there  is  no global fixed point for this  action.
Moreover,  there is  no finite orbit. The circle  $\mathcal C_1$
contains the $\omega_{h_k}$-limit set of any point in
$\TT^2\setminus \mathcal C_ 2$ and  $\mathcal C_2$ contains  the
$\alpha_{h_k}$-limit set of any point in $\TT^2\setminus\mathcal C_1
$.  Hence, the unique minimal set $M$  for this  action is contained
$\mathcal C_1$.

If $k$ is minimal, $M$  coincides with $\mathcal C_1$, the set of
$f_0$-fixed points.

If $k$ is a Denjoy homeomorphism, $M$ is strictly contained in
$\mathcal C_1$, the set of $f_0$-fixed points.

\end{itemize}

\bigskip

\noindent{\bf Proof of Corollary \ref{C2}.}

Consider BS-actions generated by $f_0$ and  $h_{\epsilon}$ given by $h_{\epsilon}(x,\theta)=(nx,  \theta+\ln(n)+\epsilon)$. If $\ln(n)+\epsilon$ is rational, then  the
 restriction  of $h_{\epsilon}$ to  $ \infty \times S^1$ is of finite order and
 every point in  $\infty \times S^1$ has a finite  BS-orbit, this action is clearly not topologically conjugate to the standard one. But, this can occur with $\epsilon$  arbitrary small, so for  a BS-action arbitrary $C^{\infty}$-close to the standard action.

\subsection {Other examples of actions of $BS(1,n)$.}  \

In this part, we construct diffeomorphisms $f$, $h$ [resp. $F$ and $H$] generating a faithful $BS(1,n)$ action on the circle [resp. the torus] where $f$ [resp. $F$] has not fixed points  but it has periodic points.

\subsubsection {On the circle.}

Let us denote $<\bar f_i, \bar h_i>$ the renormalization to $[\frac{i}{n}, \frac{i+1}{n} ]$ of the standard $BS(1,n)$-action on $\mathbb R \cup \{\infty\}$, where $i \in \{0, ..., n-1\}$.

We define $\hat f  : [0,1]/{\scriptscriptstyle (0\sim 1)}   \rightarrow [0,1]/{\scriptscriptstyle (0\sim 1)}   $ by $\hat f (x) = \bar f_i(x)$, if $x\in [\frac{i}{n}, \frac{i+1}{n} ]$ and  analogously $ h  : [0,1]/{\scriptscriptstyle (0\sim 1)}   \rightarrow [0,1]/{\scriptscriptstyle (0\sim 1)}   $ by $ h(x) = \bar h_i(x)$, if $x\in [\frac{i}{n}, \frac{i+1}{n} ]$. It is easy to see that the group generated by $\hat f$ and $h$ is isomorphic to $BS(1,n)$.

\medskip

Let $f = R_{\frac {1}{n-1}}\circ \hat f$, where $R_{\frac {1}{n-1}} (x) = x + \frac {1}{n-1} (mod 1)$.

We claim that the group generated by $f$ and $h$ is isomorphic to $BS(1,n)$.

More precisely,  $h\circ f \circ h^{-1} = h\circ R_{\frac {1}{n-1}}\circ \hat f \circ h^{-1} =  R_{\frac {1}{n-1}}\circ h \circ \hat f \circ h^{-1} $ since by construction $ R_{\frac {1}{n-1}}$ commutes with  $h$ (and also with  $\hat f$).

Then  $h\circ f \circ h^{-1} =  R_{\frac {1}{n-1}}\circ \hat f ^n =  (R_{\frac {1}{n-1}}\circ \hat f )^n = f^n$  since $R_{\frac {1}{n-1}}$ commutes with  $\hat f$  and has order $n-1$. Hence,  $f$ and $h$ generate  an  action of  $BS(1,n)$.

\smallskip

This  action is faithful, since it is a well known fact that for a non faithful action of $BS(1,n)$,  $f$ has finite order. By construction $f$ admits exactly $n-1$ periodic points of period $n-1$, so $f$ is not of finite order.

\smallskip

This construction provides an example of two circle diffeomorphisms $f$ and $h$ generating  a faithful action of  $BS(1,n)$, where $f$ has no fixed points but periodic ones.

\subsubsection {On the Torus.}

Let $f$ and $h$ be the circle diffeomorphisms as below. We define two torus diffeomorphisms :

$$\displaystyle F:\left\{ \begin{array}{ll}  (\mathbb R \cup \{\infty\} )\times [0,1]/{\scriptscriptstyle (0\sim 1)}  &\rightarrow (\mathbb R \cup \{\infty\})\times [0,1]/{\scriptscriptstyle (0\sim 1)}  \\                                                                            (x,y) & \mapsto (x + 1  , f(y) ) \end{array} \right. $$ and
$$\displaystyle H:\left\{ \begin{array}{ll}  (\mathbb R \cup \{\infty\} )\times [0,1]/{\scriptscriptstyle (0\sim 1)}  &\rightarrow (\mathbb R \cup \{\infty\})\times [0,1]/{\scriptscriptstyle (0\sim 1)}  \\                                                                            (x,y) & \mapsto (nx  , h(y) ) \end{array} \right.  . $$

The diffeomorphisms $F$ and $H$ generate a faithful action of $BS(1,n)$ on the torus, $F$ admits periodic points but not fixed points.

\section{ Isotopy class of $f$ and rotation set.}

The aim of this  section is proving Theorem \ref{rot}.
\medskip

 \begin{prop}\label{isoto} Let  $<f,h>$ be a faithful representation of $BS(1,n)$ on
   $\TT^2$. There exists a positive integer $N$ ($N\in \{1,2,3,4,6\}$)  such
   that  $f^N$ is isotopic to   identity.
\end{prop}

\begin{proof}

For proving the proposition,  it is enough to prove that there exists $N\in \mathbb{N}$ such that $A_f ^N = Id$.

\medskip

As  $A_f \in GL(2,\ZZ)$ and $f$ is conjugated to $f^n$ we have :

\begin{itemize}

\item the linear maps $A_f $ and  $A_{f ^n} = A_f ^n$ are conjugated by $A_h\in
GL(2,\ZZ)$,

\item  the  modulus of the eigenvalues of  $A_f $  are  $1$.

 \item the product of the eigenvalues is $+1$ or $-1$.

\item the trace of  $A_f $  is an integer.
\end{itemize}

\medskip

\noindent {\bf Case 1: $A_f$ admits a real eigenvalue. }

In this case, the possible eigenvalues are $+1$ or $-1$ and $A_f$ is
conjugated to one of the following applications:

\smallskip

 \hfil $A_1= \left( \begin{array}{cc}  \varepsilon_1 &  0 \\  0  &
    \varepsilon_2 \end{array}\right)$ \hfil or \hfil $A_2 = \left( \begin{array}{cc}  \varepsilon_1 &  1 \\  0  &
    \varepsilon_1 \end{array}\right)$

\noindent where $\varepsilon_i \in \{-1, 1\}$, $i=1,2$.

\medskip

It is clear that  $A_ 1 ^2 = Id$. We are going to prove that $A_2$  cannot occur.

\medskip

If $\varepsilon_1 =1$ then $A_2 ^n = \left( \begin{array}{cc}  1 & n
\\ 0 &
    1  \end{array}\right)$. One can see that  $A_2 ^n$ can not be
conjugated to $A_2$ in $GL(2,\ZZ)$. More precisely, one can compute  the
conjugating matrix  in $GL(2,\RR)$,  it is of the form  : $A_h ^{-1}= \left( \begin{array}{cc}  \frac{1}{\sqrt{n}} &  b \\  0  &
   \sqrt{n}\end{array}\right)$ where $b\in \ZZ$. This matrix does not belong
to  $GL(2,\ZZ)$. This case is not possible.

\medskip

If $\varepsilon_1 =-1$ then $A_2 ^n= \left( \begin{array}{cc}  1 & -
n   \\ 0 &
    1  \end{array}\right)$, if $n$ is even or  $A_2 ^n =
\left( \begin{array}{cc}  -1 &    n  \\ 0 &
    -1  \end{array}\right)$, if $n$ is odd.

For $n$ even,  one can easily see that  $A_2 ^n$ can not be
conjugated to $A_2$, since $Tr( A_2 ^n) = 2 \not= Tr( A_2 )$.

For $n$ odd,  one can see that  $A_2 ^n$ can not be
conjugated to $A_2$ in $GL(2,\ZZ)$:  one  compute  the
conjugating matrix in $GL(2,\RR)$,  it is of the form:
 $A_h= \left( \begin{array}{cc} \sqrt{n}  &  b \\  0  &
  \frac{1}{\sqrt{n}} \end{array}\right)$
where $b\in \ZZ$. This matrix does not belong
to  $GL(2,\ZZ)$. This case is not possible.

\medskip

\bigskip

\noindent {\bf Case 2: $A_f$ has complex eigenvalues.}

Necessary  $A_f$  has two eigenvalues
$\lambda$, $\bar \lambda$. Moreover  $\vert \lambda \vert =1 $. So $det A_ f =
\lambda\bar\lambda =  \vert \lambda \vert ^2  =1 $.

Hence, $A_f$ is conjugated to a rotation of angle $\theta$. The trace of
$A_f$ is $2\cos\theta$ and it is an integer. Then the possible values for $\cos
\theta$ are : $0, 1, -1, \frac{1}{2},-  \frac{1}{2}$.

If $\cos \theta \in \{ 1, -1\}$   then $A_f  = \left( \begin{array}{cc}  \varepsilon &  0 \\  0  &
    \varepsilon\end{array}\right)$, where $\varepsilon \in \{-1, 1\}$.

If $\cos \theta = 0$  then $A_f $ is conjugated to $ \left( \begin{array}{cc}
    0& \varepsilon_1    \\      -\varepsilon_1 & 0 \end{array}\right)$
 (where $\varepsilon_1 \in \{-1, 1\}$) which
is of order $4$, so $A_f ^4 = Id$.

If $\cos \theta \in \{ \frac{1}{2},-  \frac{1}{2} \} $ then
$A_f$ is conjugated to the  rotation of angle $ l \frac {\pi}{3}$, where $l\in
\{ 1, 2, 4 , 5\}$.  Therefore,   $ A_f ^6 = Id$.

\end{proof}

\bigskip


According to the previous proposition, given an action of
$BS(1,n) =<f,h>$ on $\TT ^2$ there
exists an integer $N$ such that $f^N$ is isotopic to identity. From now on,
we assume  that $f$ is isotopic to identity : this  is not a restrictive
hypothesis since  the action of $<f^N,h>$ on $\TT ^2$  is an action of the
Baumslag-Solitar group  $BS(1,nN)$.

\smallskip

For proving Theorem  \ref{rot}, we begin by proving the following:

\begin{prop} \label{iso}
If $f$ is isotopic to identity  and $\tilde f$ is a lift of $f$, then $\rho
(\tilde f) $ is a rational point.
\end{prop}

For proving this proposition  we need the following lemmas:

\begin{lemm}
Let $H \in Homeo_{\mathbb Z ^2}(\RR ^2)$ and $F \in Homeo^0_{\mathbb Z ^2}(\RR ^2)$  then $\rho (H \circ F \circ
H^{-1}) = A_H (\rho (F))$.
\end{lemm}

{\bf Proof}

\noindent {\bf Case 1: $A_H = Id$.} We prove that  $\rho (H \circ F \circ H^{-1}) = \rho (F)$.

Let  $(a,b)$ be a vector in the
rotation set of $H\circ F \circ
H^{-1}$. By definition,

$$(a,b) = \lim_{i\to \infty} \frac {( H\circ F \circ
H^{-1}) ^{n_i} (\widetilde{x_i})- \widetilde{x_i}} {n_i}. $$

\noindent Then $\displaystyle (a,b)=\lim_{i\to \infty}\frac{ H\circ
F ^{n_i} \circ H^{-1}(\widetilde{x_i})-\widetilde{x_i} } {n_i} = $

\smallskip

\noindent $\displaystyle=\lim_{i\to \infty} \frac { H\circ F ^{n_i}
\circ H^{-1}(\widetilde{x_i})  - F ^{n_i}\circ H^{-1}
(\widetilde{x_i})} {n_i} + \frac { F ^{n_i} \circ
H^{-1}(\widetilde{x_i})-\widetilde{x_i}} {n_i}=$

$\displaystyle \ \ \ =  \lim_{i\to \infty} \frac { (H -Id) ( F
^{n_i} \circ H^{-1}(\widetilde{x_i})) } {n_i} + \frac { F ^{n_i}
\circ H^{-1}(\widetilde{x_i})-\widetilde{x_i}} {n_i}.$

\bigskip

As $A_H = Id$, the map $(H-Id)$ is bounded (periodic) so the limit:
$$ \lim_{i\to \infty} \frac { (H -Id) ( F ^{n_i} \circ
H^{-1}(\widetilde{x_i})) } {n_i} = (0,0).$$

Moreover,

$$ \lim_{i\to \infty} \frac { F ^{n_i} \circ
H^{-1}(\widetilde{x_i})-\widetilde{x_i}} {n_i} = \lim_{i\to \infty}
\frac { F ^{n_i} \circ
H^{-1}(\widetilde{x_i})-H^{-1}(\widetilde{x_i})} {n_i}+ \frac
{H^{-1}(\widetilde{x_i})-\widetilde{x_i}} {n_i}.$$

 By definition of the rotation set, the limit :
 $$ \lim_{i\to \infty} \frac { F ^{n_i} \circ
H^{-1}(\widetilde{x_i})-H^{-1}(\widetilde{x_i})} {n_i} $$ belongs to
$\rho (F)$.

As $A_{H^{-1}} = Id$, the map $(H^{-1}-Id)$ is bounded  so the limit
:
$$ \lim_{i\to \infty} \frac { (H^{-1} -Id)
(\widetilde{x_i}) } {n_i} = (0,0).$$

 Finally, $(a,b) \in \rho
(F)$. This proves the inclusion $\rho(H\circ F \circ H^{-1}) \subset
\rho (F)$. By writing this inclusion with $H^{-1}$ instead of $H$
and $H\circ F \circ H^{-1}$ instead of $F$  we obtain  :
$\rho(H^{-1}\circ  (H\circ F \circ H^{-1})\circ H )\subset \rho
(H\circ F \circ H^{-1})$ that is  $ \rho (F) \subset \rho(H\circ F
\circ H^{-1})$.

\bigskip

\noindent {\bf Case 2:}  $H$ is a linear map that is $H =A_H $. We prove that  $\rho (H \circ F \circ H^{-1}) =A_H( \rho (F))$.

Let  $(a,b)$ be a vector in the
rotation set of $H\circ F \circ
H^{-1}$. By definition,

\medskip

$$\displaystyle (a,b) = \lim_{i\to \infty} \frac { H\circ F ^{n_i}
\circ H^{-1}(\widetilde{x_i}) -\widetilde{x_i}} {n_i} =  \lim_{i\to
\infty}   \frac { A_H\circ F ^{n_i} \circ
A_{H}^{-1}(\widetilde{x_i}) -\widetilde{x_i}} {n_i}=$$ $$=\lim_{i\to
\infty}   \frac { A_H \left( F ^{n_i} \circ
A_{H}^{-1}(\widetilde{x_i}) - A_{H}^{-1}(\widetilde{x_i})\right)}
{n_i}= A_H \left(\lim_{i\to \infty}   \frac {  F ^{n_i} \circ
A_{H}^{-1}(\widetilde{x_i}) - A_{H}^{-1}(\widetilde{x_i})}
{n_i}\right)\in A_H (\rho (F)).$$

\medskip

 This proves the inclusion $\rho(H\circ F \circ
H^{-1}) \subset A_H(\rho (F))$.  We obtain  the other
inclusion with analogous arguments as in  case 1.

\bigskip

\noindent {\bf Case 3:}  General case.

We claim that the map $A_H^{-1}\circ H \in Homeo^0_{\mathbb Z^2}(\RR ^2)$:

Let $P$ an  integer vector in $\RR ^2$ and $x$ be a point of $\RR ^2$.

$A_H ^{-1} \circ H ( x + P) = A_H ^{-1} ( H ( x)  + A_H (P))=  A_H ^{-1} \circ
H ( x)  +  A_H ^{-1} \circ A_H (P) = A_H ^{-1} \circ
H ( x)  +  P$.

By case 1, we have $\rho (A_H ^{-1} \circ H \circ F  \circ H ^{-1}  \circ A_H
) = \rho (F)$.

By case 2, we have $\rho (A_H ^{-1} \circ H \circ F  \circ H ^{-1}  \circ A_H
) = A_H ^{-1}(\rho  (H \circ F  \circ H ^{-1}) )$.

Then  $\rho (F) = A_H ^{-1}(\rho  (H \circ F  \circ H ^{-1})) $, it follows that
$\rho  (H \circ F  \circ H ^{-1}) =  A_H (\rho (F))$. \hfill $\square$

\begin{lemm} \label{trans}
 Let  $<f,h>$ be a faithful representation of $BS(1,n)$ on
   $\TT^2$ such that $f$ is isotopic to identity.  Then $\rho (\tilde f) =
   \frac {1}{n} (\tau _ Q \circ A_ {\tilde h} )(\rho (\tilde f)) $, where $Q$ is
   an integer vector in $\RR ^2$ and $\tau _Q$ denotes the translation of
   vector $Q$.

\end{lemm}

\begin{proof}

As two lifts of a  torus map differ by  an integer vector, we  have that

$\widetilde{h\circ f\circ h^ {-1}  }    = \widetilde{h} \circ
\widetilde{f}   \circ     \widetilde{ h}^ {-1}   +P $, for some integer vector $P$. By iterating this formula we
have:

$\widetilde{h\circ f^k\circ h^ {-1}} = \widetilde{h} \circ
\widetilde{f}^k\circ\widetilde{ h}^ {-1} + k P $.  Then $$\frac {\widetilde{h\circ f^k\circ h^ {-1}}}{k} = \frac {\widetilde{h} \circ
{\widetilde{f} } ^k\circ\widetilde{h}^ {-1}  }  {k}  +  P $$

Hence, by properties of the rotation set we have  $\rho (\widetilde{h\circ f\circ h^ {-1}  }   ) = \rho  (\widetilde{h} \circ
\widetilde{f}   \circ     \widetilde{ h}^ {-1}) + P = A_{\tilde h }( \rho
(\tilde f)) + P$, because of the previous lemma.

Since $f ^n = h\circ f \circ h ^{-1}$, we have  $\rho ( \widetilde {f^n}) =
A_{\tilde h }( \rho (\tilde f)) + P$.

Since $\widetilde {(f^n)} $ and $(\tilde f )^n$ are two lifts of $f^n$,  then $n \rho ( \widetilde {f}) + P' = A_{\tilde h }( \rho
(\tilde f)) + P$, for some integer vector $P'$.

Finally,   $ \rho ( \widetilde {f})  = \frac {1} {n} (\tau _
{Q}\circ A_{\tilde h })( \rho (\tilde f))$,  for some integer vector
$Q$.
\end{proof}

\noindent {\bf Proof of the Proposition \ref{iso}.}

Let $B$ be the affine map of $\RR ^2$ given by $B =  \frac {1} {n}(
\tau _ {Q} \circ A_{\tilde h })$. Note that the linear part $\bar B$
of $B$ satisfies $det \bar B = \frac {1}{n^2} det  A_{\tilde h } =
\pm \frac{1}{n ^2}$.
 The  formula given by  Lemma \ref{trans} can be written as $ \rho ( \widetilde
{f})  = B(\rho ( \widetilde {f}))$. By taking the volumes, we get:
 $ vol(\rho ( \widetilde
{f}) )  = \vert det \bar B \vert vol(\rho ( \widetilde {f})) =
\frac{1}{n^2 } vol(\rho ( \widetilde {f}))$. Then $vol (\rho (
\widetilde {f}) )  = 0$ since $\rho ( \widetilde{f}) $ is a compact
set.

This implies that   $\rho ( \widetilde{f}) $ has empty interior, so since it is a convex set, it is
either a segment or a point.

In the case where  $\rho (\widetilde{f}) $ is a segment $ [a,c]$ with $a\not=
c$, since $B([a,c]) =[a,c]$ we either have   $B(a) =a$ and  $B(c) =c$ or
 $B(a) =c$ and  $B(c) =a$ in both cases   $B^2(a) =a$ and  $B^2(c) =c$.

As $B^2$ is an  affine map, its linear part has $1$ as eigenvalue,
its trace is  $  \frac {1} {n^2}Tr (A_{\tilde h }^2 ) $ so it has
the form $\frac {p} {n^2}$ with $p\in \mathbb Z$.

Its determinant is $\frac {1} {n^4}$ so its other eigenvalue is $
\frac {1} {n^4}$.

Therefore its trace is $ 1+ \frac {1} {n^4}$ so has not the form
$\frac {p} {n^2}$ with $p\in \mathbb Z$, this is a contradiction.

Consequently, the rotation set   $\rho ( \widetilde{f}) $  is a single point
  which is the unique fixed point of the  affine  map $B$. Since $B$  has
  rational coefficients, then   $\rho ( \widetilde{f}) $ has rational
coordinates.

 \hfill $\square$

\medskip

\noindent {\bf Proof of the Theorem \ref{rot}.}

According to Proposition \ref{isoto}, there is an integer $N$ such
that  $f^N$ is isotopic to identity. By Proposition \ref{iso}, the
rotation number of any lift $\widetilde{f}^N$ is a rational vector.

  Let us write $\rho ( \widetilde{f}^N)= (\frac {p_1}{q}, \frac
{p_2}{q})$, where $p_1,p_2,q$ are integers.  Hence,  $\rho (
\widetilde{f}^{Nq})= ( {p_1}, {p_2}) \in \ZZ ^2$, then there is a
lift of $f^{Nq}$  which has rotation set equal to $\{(0,0)\}$.

According to Corollary 3.5 of \cite{MZ89},  $\{(0,0)\}= \rho (
\widetilde{f^{Nq}})=Conv( \rho_{erg} ( \widetilde{f^{Nq}}))$, where
$\rho_{erg} ( \widetilde{f}):= \{ \int (\widetilde{f} -id) d\mu,
{\text {  where } }\mu {\text { is an ergodic $f$-invariant measure}
}\}$.

Hence  $\{(0,0)\}= \rho_{erg} ( \widetilde{f^{Nq}})$. Then, using
Theorem 3.5 of \cite{Fr},  $f^{Nq} $ has a fixed point and therefore
$fix (f^{Nq})$ is non empty.

\hfill $\square$


\section{Existence of a BS-minimal set in $per(f)$.}

The aim of this section is to show the existence of a minimal set
for the action included in the set of $f$-periodic points. In the case
that $<f,h>$ is a representation of $BS(1,n)$ in $Homeo(X)$, where
$X$ is a compact metric space (Theorem \ref{minBS}), we ask for the
existence of fixed or periodic points of $f$ and in Theorem \ref{mu}
we assume that any $f$-invariant probability measure has   support
included in the set of $f$-fixed points. In this case, we also study
$f$-minimal sets and the topological entropy of $f$. In this section
we also prove Corollaries \ref{nomin} and  \ref{C1}.
\medskip

\noindent {\bf We are going to prove Theorem \ref{minBS}.}

\medskip

\noindent {\bf Proof of }(a)(1).

 Let $x_0\in fix(f)$. Since  $h^{j}\circ f\circ
h^{-j}(x_0) =f^{n^j}(x_0) = x_0$   for any $j \in \NN$,  then  $ f\circ
h^{-j}(x_0) = h^{-j}(x_0)$ so       $h^{-j}(x_0)\in  fix(f)$
  and the $\alpha$-limit set of $x_0$ for $h$ is included in $fix (f)$.

\medskip

\noindent {\bf Proof of }(a)(2).

Let $M$ be an $f$-invariant set, then  $f
(h^{-1}(M))=h^{-1}(f^n(M)) \subset h^{-1}(M)$. Since  $h^{j}\circ f\circ
h^{-j}=f^{n^j}$ for any $j\in \NN$, then $h^{-j}(M)$ is $f$-invariant.

 Let $P= fix(f)$.  It holds that
 $h^{-1}(P) \subset P$  so $h^{-j}(P)$ is a closed $f$-invariant set for any $j\in \NN$.

Let $\displaystyle K_n=\bigcap_{-n}^{n} h^{-l}(P)$, since $fix(f)\neq \emptyset$, then  $\{K_n\}$ is a family of decreasing
$f$-invariant non-empty closed sets, therefore $\displaystyle K=\bigcap_{-\infty}^{\infty}
h^{-l}(P)$ is a closed non-empty set invariant by $f$ and $h$. As a
consequence there exists a $BS$-minimal set included in $K$.

 Let us call $M_{BS}\subset fix(f)$ a minimal set for the group. Since $M_{BS}$
is $h$-invariant  there exists an $h$ minimal set $M_h \subset M_{BS}$.  The
set $M_h$ is  also $f$-invariant (it is contained in $fix (f)$), so it  follows that $M_h=M_{BS}$.

\medskip

\noindent {\bf Proof of }(a)(3).

Since $fix(f) \subset fix(f^n)$ and $h(fix(f))=fix(f^n)$ then
$\sharp\{fix(f)\}=\sharp\{fix(f^n)\}$. Therefore $fix(f)=h(fix(f))$ and $h$
has a periodic point in $fix(f)$.

\medskip

\noindent {\bf Proof of }(a)(4).

 Let $\mathcal{M}$ be a $BS$-minimal set verifying
$\mathcal{M} \cap fix(f)\neq \emptyset$. Let $ x \in
\mathcal{M} \cap fix(f)$ then $\alpha_h(x)$, the  $\alpha$-limit set of $x$ for $h$,  is a  $h$-invariant closed set  verifying   $\alpha_h(x) \subset fix(f)$. Let $M_x\subset \alpha_h(x)$ be a minimal set for the group. Since $x \in \mathcal{M}$ then $\alpha_h(x) \subset \mathcal{M}$, therefore $M_x=\mathcal{M}$ and the claim follows.

\medskip

\noindent {\bf Proof of }(b).

Suppose that $per(f) \not= \emptyset$, then there exists an positive
integer $N$ such that $fix f^N\not=  \emptyset$.

According to item (a)(2), there is a minimal set $M_N$ of $<f^N, h>$
such that $M_N \subset fix f^N$.


Let $\mathcal M = \cup_{k=0} ^{N-1} f^k(M_N)$.

\smallskip

\begin{itemize}

\item  $\mathcal M$ is $f$-invariant : $\displaystyle f(\mathcal M)= \cup_{k=0} ^{N-1} f^{k+1}(M_N)= \mathcal M$, since $f^N (M_N) = M_N$.
In fact,  $ \displaystyle \mathcal  M = \cup_{k\in  \mathbb Z}
f^k(M_N)$.

\item $\mathcal  M$ is $h$-invariant :  $h(\mathcal M)= \cup_{k=0} ^{N-1} h\circ f^{k}(M_N)=  \cup_{k=0} ^{N-1} f^{nk}\circ h (M_N)=
\cup_{k=0} ^{N-1} f^{nk} (M_N) \subset \mathcal M$, since $M_N$ is
$h$-invariant.

\end{itemize}

Since $\mathcal M$ is closed, non empty, $f$ and $h$ invariant, it
contains a BS minimal set $ M$.

\begin{itemize}
\item $\mathcal  M\subset fix (f^N)$: let $x\in \mathcal  M$, there is a $k=0, ..., N-1$ and $x'\in M_N$ such that  $x=f^k (x')$. Hence  $f^N(x) =   f^{N+k}(x')= f^{k}(f^N (x'))=   f^{k}(x')=x $.
\end{itemize}

Finally, $M\subset \mathcal  M\subset fix (f^N)$.

\hfill $\square$

\medskip

\noindent {\bf Proof of Corollary  \ref{nomin}.}

\medskip

If the action were minimal, its unique minimal set would be $X$ and
would be contained in $fix f ^N$, according to item (b) of Theorem
\ref{minBS}. This implies that $X= fix f^N$ and  so $f^N =Id$,  the
action would not be faithful. This is a contradiction. \hfill
$\square$

\bigskip

%
%
%
%
%

\noindent {\bf The following is the  proof of Theorem \ref{mu}.}

\noindent {\bf Proof of }(1). Let $M_f$ be an $f$-minimal set and
$x_0 \in M_f$. Let $\displaystyle \mu_k=
\frac{1}{k}\sum_{i=0}^{k-1}\delta_{f^i(x_0)}$  and $\mu$ a weak
limit of $\mu_k$. It is known that $\mu$ is an $f$-invariant
probability measure and  its support is included in
$M_f=\overline{\mathcal{O}_f (x_0)}$, the closure of the $f$-orbit
of $x_0$. In addition, by hypotheses $supp(\mu) \subseteq fix(f)$
then $M_f \cap fix(f) \neq \emptyset.$ It follows that $M_f$ is
reduced to a fixed point. Since a periodic orbit is a minimal set we
have that  $per(f)$ coincides with the set $fix(f)$.

\medskip

 \noindent {\bf Proof of }(2).
  Let $M_{BS}$ be an $BS$-minimal set.
 Since $M_{BS}$ is  $f$-invariant, there exists an $f$ minimal set $M_f \subset M_{BS}$. Since $M_f$ is an $f$-fixed point,
   it follows that $M_{BS}\cap fix(f)\neq \emptyset$, so according to items (a)(2) and (a)(4) of Theorem \ref{minBS}, $M_{BS} \subseteq fix(f)$
   and it coincides with a $h$-minimal set in $fix(f)$.
\medskip

 \noindent {\bf Proof of }(3).
 Recall that $ent_{top}(f)= sup\{ ent_{\nu }(f)\}$ where  $\nu$ is an $f$-{invariant  probability measure}, the supremum of all metric entropies.
 Since $ ent_{\nu }(f)=ent_{\nu }(f \vert_{supp(\nu)})$ and $f \vert_{supp(\nu)}=Id$, we have that $ ent_{top}(f)=0$.

\hfill $\square$

\bigskip
\noindent {\bf We finish this section by proving Corollary \ref{C1}.}

\medskip

If $S=\TT^2$, let $N$ be the positive integer given by Theorem 1, then
the set $fix(f^N)\neq \emptyset$. If $S=S^2$, let $N$ be the
smallest positive integer such that $ f^N$ has at least three fixed
points and it is orientation preserving. Otherwise, $N=1$.

 In addition, $f$ is distortion element of $BS(1,n)$, so according to Theorem 1.3 of \cite{FH} for any $f$-invariant probability measure, $\mu$, it holds that  $supp(\mu) \subseteq fix(f^N)$ so Theorem \ref{mu} implies the claim of this corollary except the $f$-ellipticity of the points in minimal sets.

 For simplicity, we will prove the ellipticity in the case where $f$ has fixed points. The general case is analogous.

 Let $x_0$ be a point in a $BS$-minimal set $M_{BS}$. Since  $M_{BS}$ is also an $h$-minimal set, the $h$-orbit of $x_0$ is recurrent. Then there exists a subsequence $(n_k)$ $(n_k \rightarrow \infty)$ such that $h ^{-n_k } (x_0) \rightarrow x_0$.

From  $h  ^{n_k }\circ f \circ h^{-n_k } = f ^{n^{n_k }}$, we deduce
that :$$Dh ^{n_k } (f(h^{-n_k }(x_0)))  \circ Df (h^{-n_k
}(x_0))\circ Dh^{-n_k }(x_0)  = Df^{n^{n_k }}(x_0).$$

As the points $x_0$ and $h^{-n_k }(x_0)$ are fixed by $f$ and
$(Dh^{-n_k }(x_0)) ^{-1} =Dh ^{n_k } (h^{-n_k }(x_0)) $, then:
$$(Dh^{-n_k }(x_0)) ^{-1} \circ Df (h^{-n_k }(x_0))\circ Dh^{-n_k }(x_0)  = (Df(x_0))^{n^{n_k }}.$$

So  $Df (h^{-n_k }(x_0))$ and $(Df(x_0))^{n^{n_k }}$ have the same
eigenvalues. As $f$ is $C^ 1$ and  $h ^{-n_k } (x_0) \rightarrow
x_0$, then  $Df (h^{-n_k }(x_0))\rightarrow Df(x_0)$.

We conclude that   $Df(x_0)$ and $(Df(x_0))^{n^{n_k }}$ have the
same eigenvalues,  finally the eigenvalues of $Df(x_0)$ have module
$1$.

\hfill $\square$

\section {Perturbations of the standard $BS(1,n)-$action  on $\TT^2$.}

\noindent  Let us recall that :

 $\bullet$ the standard BS-action on $\TT ^2$ is the one  generated by the two diffeomorphisms of $\mathbb R \cup \{\infty\} \times S^1$ :
 $f_0(x,\theta)=(x+1, \theta)$ and $ h_0(x,\theta)= (nx,\theta+\ln(n)), $

$\bullet$  $\mathcal C_1:= \infty \times S^1$ and  $\mathcal C_2:= 0
\times S^1$.

\medskip

\bigskip
Before proving  Theorem \ref{minset}, we prove the following
\begin{lemm}
\label{pers}
Let us consider a BS-action $<f,h>$ on $\TT ^2$ generated by  $f$
and $h$ sufficiently $C^0$-close to homeomorphisms $\bar f_0$ and $\bar h_0$ which generate a  BS-action.
If both  $\bar f_0$ and $\bar h_0$ are isotopic to identity and the rotation set  of a lift of $\bar f_0$ is $(0,0)$, then  the rotation set of a lift of $f$ is $(0,0)$.

\end{lemm}

\smallskip

\noindent {\bf Proof of the lemma.}

For $(f,h)$  sufficiently close to $(\bar f_0, \bar h_0)$, $f$ and $h$ are isotopic to identity. By Lemma \ref{trans} the rotation set of  any  lift $\tilde f$ of $f$ satisfy    $n \rho (\tilde f) = \rho (\tilde f) + (p,q) $, where $(p,q)$ is  an integer vector.  Then  the rotation set of $\tilde f$ is a rational vector $(\frac {p}{n-1}, \frac {q}{n-1})$, with $p, q$ integers.

It is proved in  \cite{MZ89} that  the rotation set map $\rho : Homeo_{\mathbb Z ^2}(\RR ^2) \rightarrow \mathcal K(\RR^2)$, the set of  compact subsets of $\RR^2$  is upper semi-continuous with respect  to the compact-open topology on $Homeo_{\mathbb Z ^2}(\RR ^2)$ and the Haussdorff topology on $\mathcal K(\RR^2)$. In other words, if $G$ is an element of $Homeo_{\mathbb Z ^2}(\RR ^2)$ and $U$  is a neighborhood of $\rho(G)$ in $\RR^2$, then  for $F$ sufficiently close to $G$, we have $\rho(F) \subset U$.

The rotation set of a lift  of $\bar f_0$ is $(0,0)$, consider a neighborhood $U$ of $(0,0)$ in $\RR^2$ that  contains no points of the form $(\frac {p}{n-1}, \frac {q}{n-1})$, with $p, q$ integers and $(p,q)\not=(0,0)$.

According to the previous result of \cite{MZ89},  for $f$ sufficiently close to $f_0$, the rotation set of a lift of $f$ is included in  $U$ and  since it has  form $(\frac {p}{n-1}, \frac {q}{n-1})$, it must be $(0,0)$. \hfill $\square$

\bigskip

\noindent {\bf Proof of  Theorem \ref{minset}.}

\medskip

(1) The circles $\mathcal C_1$ and  $\mathcal C_2$ are $h_0$-normally
hyperbolic in the sense of \cite{HPS}. Consider a neighborhood $U_1$ of
$\mathcal C_1$  where $\mathcal C_1$ is $h_0$-attractive and  a neighborhood
$U_2$ of $\mathcal C_2$  where $\mathcal C_2$ is $h_0$-repulsive. Obviously,
there exists some integer $k_0$ such that $ h_0^{k} (\TT ^2 \setminus U_2)
\subset U_1$,  $ h_0^{-k} (\TT ^2 \setminus U_1)
\subset U_2$, for all $k\geq k_0$.

According to  Theorem 4.1 of \cite{HPS}, there exists a $C^1$-neighborhood
$\mathcal V$ of $h_0$ in $Diff^1 (\TT ^2)$ such that for all $h\in \mathcal V$
 there exist two circles  $\mathcal C'_1$ and  $\mathcal C'_2$ which are
 $C^1$-closed to $\mathcal C_1$ and  $\mathcal C_2$ respectively and they are $h$-invariant.

Moreover
 $\mathcal C'_1$ is $h$-attractive in $U_1$  and    $\mathcal C'_2$  is
 $h$-repulsive in $U_2$ and   $ h^{k} (\TT ^2 \setminus U_2)
\subset U_1$,  $ h^{-k} (\TT ^2 \setminus U_1) \subset U_2$,
 for all $k\geq k_0$. Therefore,  item (1) is proved.

\bigskip

(2)   Obviously, the rotation set of a lift of $f_0$ is $(0,0)$. According to previous lemma, the rotation set of a lift of $f$  must be $(0,0)$. Then $fix(f)$ is not empty.

 Let $x_0\in fix(f)$, if $x_0\notin \mathcal C'_1$ the its $\alpha$-limit set for
$h$ is included in $\mathcal C'_2$ and consists of $f$-fixed points, according
to Theorem \ref{minBS}. In other
words,  $\mathcal C'_2$ intersects  $fix(f)$. But for $f$ sufficiently close to
$f_0$  we have that $f^j(\mathcal C'_2) \cap \mathcal C'_2
=\emptyset$, for any $j\not=0$.  Hence $x_0\in \mathcal C'_1$.

\bigskip

(3)  We first prove that any minimal set of BS intersects $\mathcal C'_1 $.

Let us consider $M$ a BS-minimal set :

Suppose that $M \subset  \mathcal C'_2$ then $ f(M)=M \subset f(\mathcal
C'_2)$ then   $f(\mathcal C'_2 )  \cap \mathcal C'_2 \not= \emptyset$ which is
contradiction for $f$ close to $f_0$.

Since $M  \not\subset  \mathcal C'_2$,  there is
$x_0\in M \setminus \mathcal C'_2 $. Then  $\omega_h(x_0)\subset  \mathcal
C'_1\cap  M  $, so  we are done.

\medskip

 The circle $\mathcal C'_1$ is $h$-invariant, we can consider the rotation
number $\rho$ of the restriction of $h$ to  $\mathcal C'_1$ :

\medskip

\noindent {\bf Case 1:  $\mathbf \rho \in \QQ$.}

\smallskip
There is a BS-minimal set $\mathcal M$ included in $fix(f)$ so in  $\mathcal C'_1$.  This  set $\mathcal M$ contains an $h$-minimal set in  $\mathcal C'_1$.   Moreover $h$ has periodic orbit and any minimal set of  $h_{\vert \mathcal C'_1 }$ is an  $h$-periodic orbit.  Then there is an $h$-periodic orbit contained in $\mathcal M \subset fix(f)$.  So this $h$-periodic orbit is a finite  BS-orbit.

\noindent {\bf Case 2:  $\mathbf \rho \notin \QQ$.}

\smallskip

{\bf Case 2a:   $h_{\vert \mathcal C'_1}$ is conjugated  to an irrational
  rotation.}

\smallskip

We claim  that  $\mathcal C'_1 = fix(f)$.  Let  $x_0 \in fix( f) $, then
$\alpha_h(x_0)= \mathcal C'_1$ and it is contained in $fix(f)$. Hence $fix(f)=\mathcal C'_1$.

\medskip

Now, we prove that $\mathcal C'_1$ is a minimal set for the BS-action:

Let $x$ be in $\mathcal C'_1 = fix f$. The closure of $h$-orbit of $x$ is BS-invariant and coincide with $\mathcal C'_1$.
Consequently,  the circle   $\mathcal C'_1$ is a minimal set for the BS-action.


\medskip

{\bf Case 2b:   $h_{\vert \mathcal C'_1}$ is semi-conjugated (not conjugated)
  to an irrational  rotation. }

\smallskip

 Then $h_{\vert \mathcal C'_1}$ admits a unique minimal set $K$ that is
 homeomorphic to a Cantor set.

Let $x_0$ be a fixed  point of $f$, then  $K=\alpha_h(x_0)\subset fix(f)$. So $K$ is BS-invariant, so it contains a BS-minimal set.


Since any  BS-minimal set intersects $\mathcal C'_1$ and $\alpha_h(x) = K$ for
all $x\in \mathcal C'_1$, any  BS-minimal set contains $K$.

Finally, $K$ is the unique  BS-minimal set and it is a Cantor set.

\bigskip

 \hfill $\square$

In the case that the action is $C^2$ we have the following

\medskip

\begin{coro}
\label{Denjoy}
If the action is $C^2$ and sufficiently $C^1$-close to $<f_0
  ,h_0>$ then  either :

\begin{enumerate}

\item  $\mathcal C'_1 = fix(f)$ is the unique minimal set for the action and the
minimal sets of $f$ are its fixed points or

\item there  exists a finite BS-orbit contained in $\mathcal C'_1$.

\end{enumerate}

\end{coro}

\noindent {\bf Proof. }

\medskip

According to theorem \ref{minset} (3), either there exists a finite $BS$-orbit
in $\mathcal C'_1$ or the action has an unique minimal set $M$ which is  the
unique $h_{\vert \mathcal C'_1}$-minimal set.

In the second case,  since $h$ is $C^2$, the circle map  $h_{\vert
\mathcal C'_1}$ is $C^2$ and according to Denjoy's  theorem, $M$ is
the whole circle   $\mathcal C'_1=fix(f)$. \hfill $\square$

\medskip


\subsection{Persistent global fixed point}

 \begin{prop}

Let us consider a BS-action $<f,h>$ on $\TT ^2$ generated by  $f$ and $h$ sufficiently $C^1$-close to $\bar f_0$ and $\bar h_0$, where $\bar f_0$ and $\bar h_0$ are isotopic to identity.   If the rotation  set of a lift of $\bar f_0$  is $(0,0)$ and $\bar h_0$ is a Morse Smale diffeomorphism satisfying that any periodic point is $\bar h_0$-fixed,  then $<f,h>$ admits fixed point.

\end{prop}

\medskip
\noindent {\bf Proof.}

Any $h$ sufficiently $C^1$-close to $\bar h_0$ is a Morse Smale diffeomorphism where any $h$-periodic point is fixed. In particular, any $h$-minimal set is an $h$-fixed point.

By lemma \ref{pers}, rotation set of a lift $f$ is $(0,0)$, so $fix(f)$ is not empty.

As a consequence of Theorem \ref{minBS}, there is a BS-minimal set included in  $fix(f)$, this minimal set contains an $h$-minimal set, that is a fixed point of $h$. This point is a global fixed point.$\hfill \square$

\bigskip

\noindent {\bf Explicit example.}

Let   $\bar f_0(x,\theta)=(x+1,\theta +1)$ and $\bar h_0(x,\theta)=( nx, n\theta),$ where $x\in \RR \cup \infty$ and $\theta \in \RR \cup \infty$. It is easy to check that both diffeomorphisms are isotopic to identity, $\bar h_0$ is a Morse Smale diffeomorphism with two fixed points:  $(0,0)$ and  $(\infty,\infty)$ and  $\bar f_0$ has
a unique fixed point:   $(\infty,\infty)$. These diffeomorphisms satisfy the hypothesis of the previous proposition, so any sufficiently $C^1$-close BS- action has  fixed point.









\medskip

\bibliographystyle{alpha}
\bibliography{refenBS1}

\end{document}